\documentclass[10pt]{amsart}

\usepackage{amsxtra}
\usepackage{amssymb}

\newcommand\la{\langle}
\newcommand\ra{\rangle}

\newcommand\dd{{\mathfrak d}}
\renewcommand\aa{{\mathfrak a}}

\newcommand\bb{{\mathfrak b}}

\newcommand\zz{{\mathfrak z}} 
\newcommand\rr{{\mathfrak r}}
\newcommand\hh{{\mathfrak h}}
\newcommand\jj{{\mathfrak j}} 

\newcommand\nn{{\mathfrak n}} 
\newcommand\ggo{{\mathfrak g}}
\newcommand\vv{{\mathfrak v}}

\newcommand\ee{{\mathfrak e}}
\newcommand\aff{{\mathfrak{aff}}}
\newcommand\glo{{\mathfrak{gl}}}

\newcommand\CC{\mathbb C}

\newcommand\RR{\mathbb R}

\newcommand\GL{{\mathrm{GL}}}
\newcommand\End{{\mathrm{End}}}

\newcommand\ad{\operatorname{ad}}

\renewcommand\Im{\operatorname{Im}}

\theoremstyle{plain}

\newtheorem{thm}{Theorem}[section]
\newtheorem{lem}[thm]{Lemma}
\newtheorem{prop}[thm]{Proposition}
\newtheorem{cor}[thm]{Corollary}

\theoremstyle{definition}
\newtheorem{defn}[thm]{Definition}
\newtheorem{remark}[thm]{Remark}

\begin{document}
\title[symplectic Lie algebras]{four dimensional symplectic Lie algebras}

\author{Gabriela Ovando}
\thanks{G. Ovando: CIEM - Facultad de Matem\'atica, Astronom\'\i a y F\'\i sica,
Universidad Nacional de C\'or\-do\-ba, C\'or\-do\-ba~5000, Argentina \\
email: ovando@mate.uncor.edu}

\thanks{{\it (2000) Mathematics Subject Classification}: Primary:53D05; Secondary: 22E25, 17B56 }

\thanks{{\it Key words and phrases:} symplectic structures, solvable Lie algebra, cotangent extension, symplectic double extensions, cohomology}

\date{\today}

\begin{abstract}In this paper we deal with symplectic Lie algebras.  All symplectic  structures are determined for dimension four and the corresponding Lie algebras are classified up to equivalence.  Symplectic  four dimensional Lie algebras are described either as solutions of the cotangent extension problem or as symplectic double extension of $\RR^2$ by $\RR$. The difference in the choice of a certain model lies on the existence or not of a lagrangian ideal. Moreover all extensions of a two dimensional Lie algebra are determined, and so all solutions (up to equivalence) of the cotangent extension problem are
given in dimension two. By studying the adjoint representation we generalize to higher dimensions finding   obstructions to the existence of symplectic forms.  Finally, as an appendix we compute the authomorphisms of four dimensional symplectic Lie algebras and the cohomology over $\RR$ of the solvable real four dimensional Lie algebras.
\end{abstract}

\maketitle

\section{Introduction}

One origin of the study of symplectic geometry is classical mechanics. Symplectic structures have proved to be an important tool in the description and geometrization of several phenomena.
Special cases of symplectic manifolds are the K\"ahler ones, in which the symplectic form derives from a complex structure and a compatible Riemannian metric.
In this context the methods of rational homotopy theory have been applied  succesfully to symplectic geometry. These techniques were useful in attacking many geometric problems such as the construction of symplectic manifolds without K\"ahler structure among others (see \cite{O-T} for a survey  on this topic and other articles such as \cite{A} \cite{A-G} \cite{F-G} \cite{F-G-G} \cite{L} \cite{Mc} \cite{Ya} for example). Other algebraic tools were used  to attack geometric aspects in \cite{B-G} \cite{DM2} \cite{DN} \cite{LM}.

In addition to K\"ahler manifolds, one of the most frequently encountered types  of symplectic manifolds is the cotangent bundle of a differentiable manifold. Adding some algebraic data  one has another class of examples provided by Lie groups, endowed with a left invariant symplectic structure, that is, left invariant translations by elements of the Lie group are symplectomorphisms. Let $H$ denote a Lie group with Lie algebra $\hh$ and let $\hh^{\ast}$ be the dual vector space of $\hh$. Let $T^{\ast}H$ be the cotangent bundle of $H$, then $H$ is identified with the zero section in $T^{\ast}H$ and $\hh^{\ast}$ with the fiber over a neutral element of $H$. Let $\omega_0$ be the skew symmetric bilinear form defined in $\hh^{\ast} \oplus \hh$ by setting
$$\omega_0((\varphi, x), (\varphi', x')) = \varphi(x') -\varphi'(x)$$
It is easy to see that $\omega_0$ spans a left invariant closed 2-form in the canonical Lie group structure of $T^{\ast}H$ if and only if $H$ is abelian (see Remark (\ref{cano})). That leads to ask the following: {\it is there a Lie group structure on $T^{\ast}H$  in such way that the left invariant two form induced by $\omega_0$  is closed?} That is known as the {\it cotangent extension problem} \cite{By1}. More precisely, let $\hh$ be a Lie algebra and let $\hh^{\ast}$ be the dual vector space of $\hh$ and $\omega_0$ be as defined above. The problem is to find a Lie algebra structure on $\hh^{\ast} \oplus \hh$ which  satisfies the following conditions:

(c1) $0 \longrightarrow \hh^{\ast} \longrightarrow \hh^{\ast} \oplus \hh \longrightarrow \hh \longrightarrow 0$ is an exact sequence of Lie algebras, $\hh^{\ast}$ endowed with the abelian Lie algebra structure;

(c2) the left invariant 2-form spanned by $\omega_0$ is closed.

Note that $\hh^{\ast}$ is a lagrangian ideal on $\hh^{\ast} \oplus \hh$. The resulting Lie algebra is called a solution of the cotangent extension problem.
This construction of symplectic manifolds was used by Boyom \cite{By1} \cite{By2} to give models for symplectic Lie algebras.

Another construction arises by the  {\it symplectic double extension} given  by  Medina and Revoy \cite{MR} and generalized by Dardi\'e and Medina \cite{DM1}. This construction realizes a symplectic group as the reduction of another symplectic group. In this situation one has a symplectic Lie algebra $\ggo$ and an isotropic ideal $\hh$, such that $\hh^{\perp}$ is an ideal and
 the symplectic structure on $\hh^{\perp}/\hh$ is induced from that of $\ggo$.

Related to symplectic structures one has many problems. One of them is the existence and classification of such structures on any given Lie algebra and another problem is the algebraic structure of a Lie algebra (group) admitting  such a symplectic structure. The first problem  was studied for example in \cite{BG} \cite{AC}. The second problem was treated for example in \cite{By1} \cite{By2} \cite{DM1} \cite{MR}. In \cite{MR} the authors determined  symplectic four dimensional Lie algebras. In this paper  we compute  all symplectic structures on any real solvable four dimensional Lie algebras  showing in this way, all symplectic Lie algebras endowed with a  symplectic structure. Then we consider the action of the authomorphism group of the Lie algebra on the space of symplectic structures and we get the classification of symplectic Lie algebras up to this equivalence relation, completing the table in \cite{MR}. In the four dimensional case we study  the action of the adjoint representation  and we get obstructions to the existence of symplectic structures in  higher dimensions. A next goal in this paper is to reconstruct these symplectic Lie algebras in terms of the above described models. In this sense the principal result we  prove is that any symplectic Lie algebra which is either completely solvable or $\aff(\CC)$ is a solution of the cotangent extension problem. Furthermore we classify all solutions of the cotangent extension problem up to equivalence for $\hh$ of dimension two, that is, $\hh$ abelian or isomorphic to $\aff(\RR)$. For the other symplectic four dimensional Lie algebras we apply results of \cite{MR} \cite{DM1} to prove that $\ggo$ is obtained as a double extension of the  two-dimensional abelian Lie algebra by $\RR$. Essentially, the choice of a certain model lies on the existence or not of a lagrangian ideal on the symplectic Lie algebra. Lie algebras having a lagrangian ideal are modelized as solutions  of the cotangent extension problem. The other Lie algebras described as symplectic double extensions have an isotropic ideal which cannot be lagrangian.  However it is also possible to describe some symplectic Lie algebras having lagrangian ideals also as symplectic double extensions. In fact this is the case of symplectic  nilpotent Lie algebras, which can be obtained by a finite sequence of (classic) symplectic double extensions (Theorem of \cite{MR}).

The paper is organized as follows: the first section shows all  symplectic four dimensional Lie algebras, the corresponding symplectic structures and its classification. To this end we make use of the classification of real solvable four dimensional Lie algebras with the notations given in \cite{A-B-D-O}. Reading the previous results we point out the exact symplectic Lie structures, which were obtained by Campoamor \cite{Ca}. Finally we describe the adjoint representation of these symplectic Lie algebras, giving  obstructions to the existence of symplectic structures in Section \ref{last}. The second section is devoted to find models for the symplectic four dimensional Lie algebras. For each model we recall the main definitions and give examples in the class of symplectic four dimensional Lie algebras. As an appendix we give the authomorphisms of four dimensional symplectic Lie algebras and we compute explicitly the real cohomology of four dimensional solvable Lie algebras.

\

If $G$ is a simply connected Lie group then its Lie algebra will be denoted with greek letters $\ggo$ and identified as usual with the left invariant vector fields of $G$.

\section{Four dimensional symplectic Lie algebras}

\subsection{Four dimensional solvable Lie algebras}
Since we are interested on left invariant structures, our work is reduced to the Lie algebras of the corresponding Lie groups. As a  first step we  exhibit in the following proposition  the different classes of four dimensional solvable Lie algebras (see \cite{D}, \cite{Mu} or \cite{A-B-D-O}).

\begin{prop} \label{clases} Let $\ggo$ be a solvable four dimensional real Lie algebra. Then if $\ggo$ is not abelian, it is equivalent to one and only one of the Lie algebras listed below:

$$\begin{array}{ll}
{\rr\hh_3:} & {[e_1,e_2] = e_3 }  \\
{\rr\rr_3:\qquad } & {[e_1,e_2] =  e_2,\,[e_1,e_3] = e_2+ e_3 }  \\
{\rr\rr_{3,\lambda}:} & {[e_1,e_2] =e_2, [e_1,e_3]=\lambda e_3} \qquad {\lambda \in [-1,1]}\\
{\rr \rr'_{3,\gamma}:} & {[e_1, e_2] = \gamma e_2 -  e_3,[e_1, e_3] =  e_2 + \gamma e_3}  \qquad {\gamma \ge 0}\\
{\rr_2 \rr_2:} & {[e_1, e_2] = e_2,\,[e_3,e_4] = e_4}\\
{\rr_2':} & {[e_1,e_3] = e_3,\,[e_1,e_4] = e_4,\,[e_2,e_3]= e_4, \,[e_2,e_4]=-e_3} \\
{\nn_4:} & {[e_4,e_1] = e_2,\,[e_4, e_2] = e_3 } \\
{\rr_4:} & {[e_4,e_1]=e_1,\, [e_4,e_2] = e_1 + e_2, [e_4, e_3]=e_2 + e_3}\\
{\rr_{4,\mu}:} & {[e_4,e_1] = e_1,\, [e_4,e_2] = \mu e_2,\,[e_4,e_3] = e_2 + \mu e_3 } \qquad {\mu \in \RR} \\
{\rr_{4,\alpha, \beta}:}& {[e_4,e_1] = e_1,\, [e_4,e_2] = \alpha e_2,\, [e_4,e_3] = \beta e_3,}\,  \\
& \text{{\rm with }} \,-1 < \alpha \leq \beta \le 1 , \, \alpha \beta \ne 0,\,  \text{\rm{or }} \,-1 = \alpha \leq \beta\le 0 \\
{\rr'_{4,\gamma, \delta}:} & {[e_4,e_1] = e_1,\, [e_4,e_2] = \gamma e_2 - \delta e_3, \, [e_4,e_3] = \delta e_2 +\gamma e_3} \quad \gamma \in \RR, \delta > 0\\
{\dd_4:} & {[e_1,e_2]=e_3,\, [e_4,e_1] = e_1,\, [e_4,e_2] = -e_2}   \\
{\dd_{4,\lambda}:} & {[e_1,e_2]=e_3,\, [e_4,e_3] = e_3,\,[e_4, e_1]=\lambda  e_1,\,\,[e_4, e_2]=(1-\lambda) e_2} \quad  \lambda \ge \frac12\\
{\dd'_{4,\delta}:} & {[e_1,e_2]=e_3,\, [e_4, e_1]=\frac{\delta}2  e_1- e_2,[e_4,e_3] = \delta e_3,\,\,[e_4, e_2]=e_1+ \frac{\delta}2 e_2} \quad \delta \ge 0\\
{\hh_4} & {[e_1,e_2]=e_3,\, [e_4,e_3] = e_3,\,[e_4, e_1]=\frac{1}2 e_1,\,\,[e_4, e_2]= e_1 + \frac{1}2 e_2 }  \\
\end{array}
$$
\end{prop}

\begin{remark}
Observe that $\rr_2\rr_2$ is the Lie algebra  $\aff(\RR) \times \aff(\RR)$, where $\aff(\RR)$ is the Lie algebra of the Lie group  of affine motions of $\RR$, $\rr_2'$ is  the real Lie algebra underlying on the complex  Lie algebra $\aff(\CC)$, $\rr\rr_{3,-1}$ is the trivial extension of $\ee(1,1)$, the Lie algebra  corresponding to the Lie group of rigid motions of the Minkowski 2-space;
$\rr'_{3,0}$ is the trivial extension of $\ee(2)$, the Lie algebra of the  Lie group of rigid motions of $\RR^2$; $\rr\hh_3$ is the trivial extension of
the three-dimensional Heisenberg Lie algebra  denoted by $\hh_3$.
\end{remark}

A Lie algebra is called {\it unimodular} if tr($\ad_x$)=0 for all $x\in \ggo$, where tr denotes the trace of the map. The application $x  \to$tr($\ad_x$) is an homorphism of Lie algebras, thus its kernel is an ideal called the {\it unimodular kernel} of $\ggo$.
 The unimodular four-dimensional solvable Lie algebras algebras  are: $\RR^4, \; \;  \rr
\hh_3, \;\; \rr \mathfrak r_{3,-1}, \;\;  \rr\mathfrak
r'_{3,0},\;\;\mathfrak n_4, \;\; \mathfrak r_{4,-1/2},\linebreak \mathfrak r_{4,\mu,-1-\mu}$ {\scriptsize $(-1<\mu\leq -1/2)$}, $\;\; \mathfrak r'_{4,\mu,-\mu/2},\,\;\;\mathfrak d _{4},\;\; \mathfrak d_{4,0}'.$

Recall that a solvable Lie algebra is {\it completely solvable} when $\ad_x$ has real eigenvalues for all $x\in\ggo$.

\begin{remark} For an explanation concerning the Lie groups which admit a compact quotient, see for example the work of Oprea and Tralle \cite{O-T}. In particular if $G$ admits a discrete subgroup $\Gamma$ with compact quotient, then the corresponding Lie algebra is unimodular \cite{Mi}.
\end{remark}

\

\subsection{Classification of symplectic Lie algebras}
 A {\it symplectic structure} on a 2n-dimensional Lie algebra $\ggo$ is a  closed 2-form $\omega \in \Lambda^2(\ggo^{\ast})$ such that $\omega$ has maximal rank, that is, $\omega^n$ is a volume form on the corresponding Lie group. Lie algebras (groups) admitting symplectic structures are called {\it symplectic } Lie algebras (resp. Lie  groups).

It is known that if $\ggo$ is  four dimensional and symplectic then it must be solvable \cite{Ch}. However not every four dimensional solvable Lie group admits a symplectic structure. In this section we determine all left invariant symplectic structures on simply connected four dimensional Lie groups and we classify the corresponding Lie algebras, up to equivalence. This work completes the table of \cite{MR}, where all symplectic four  dimensional Lie algebras were determined. With this results we get exact symplectic Lie algebras and we point out the action of the adjoint representation on symplectic Lie algebras, characterizing four dimensional Lie algebras admitting  symplectic structures.

 Denoting by $\{e^i\}$ the dual basis on $\ggo^{\ast}$ of the basis $\{e_i\}$ on $\ggo$ (see(\ref{clases})),  the next  Proposition \ref{simplecticas}  describes symplectic structures in the four dimensional case.

\begin{prop} \label{simplecticas} Let $\ggo$ be a symplectic real Lie algebra of dimension four. Then $\ggo$ is isomorphic to one of the following Lie algebras equipped with a symplectic form as follows:

$$
\begin{array}{ll}
{\rr\hh_{3}:} & {\omega=a_{12}e^1 \wedge e^2+a_{13}e^1 \wedge e^3 + a_{14}e^1 \wedge e^4 + a_{23}e^2 \wedge e^3+a_{24}e^2 \wedge e^4} \, \\
& \qquad {a_{14}a_{23}-a_{13}a_{24}\ne 0} \\
{\rr \rr_{3,0}:}& \omega=a_{12}e^1 \wedge e^2 + a_{13}e^1 \wedge e^3 +a_{14} e^1 \wedge e^4 +a_{34} e^3 \wedge e^4,\, \\& \qquad {a_{12}a_{34}\ne 0} \\
{\rr \rr_{3,-1}:}& \omega=a_{12}e^1 \wedge e^2 + a_{13}e^1 \wedge e^3 +a_{14} e^1 \wedge e^4 +a_{23} e^2 \wedge e^3,\, \\
& \qquad a_{14} a_{23} \ne 0 \\
{\rr\rr'_{3,0}:} & \omega={a_{12}e^1 \wedge e^2 + a_{13}e^1 \wedge e^3+a_{14}e^1 \wedge e^4 + a_{23}e^2 \wedge e^3} \,\\ & \qquad {a_{14}a_{23}\ne 0} \\
{\rr_2\rr_2:} & \omega={a_{12}e^1 \wedge e^2 + a_{13}e^1 \wedge e^3 + a_{34}e^3 \wedge e^4},\\
& \qquad  \, {a_{12}a_{34}\ne 0} \\
\end{array}$$
$$\begin{array}{ll}
{\rr'_2:} & \omega={a_{12}e^1 \wedge e^2 + a_{13-24}(e^1 \wedge e^3 - e^2 \wedge e^4) +a_{14+23}(e^1 \wedge e^4+e^2 \wedge e^3) } \,\\
& \qquad {a_{14+23}^2 + a_{13-24}^2\ne 0} \\
{\nn_{4}:} & \omega={a_{12}e^1 \wedge e^2 + a_{14}e^1 \wedge e^4+a_{24}e^2 \wedge e^4 + a_{34}e^3 \wedge e^4},\, \\ & \qquad {a_{12}a_{34}\ne 0} \\
{\rr_{4,0}:} & \omega={a_{14}e^1 \wedge e^4 + a_{23}e^2 \wedge e^3+a_{24}e^2 \wedge e^4 + a_{34}e^3 \wedge e^4},\, \\
{\rr_{4,-1}:} & \omega={a_{13}e^1 \wedge e^3 + a_{14}e^1 \wedge e^4+a_{24}e^2 \wedge e^4 +a_{34}e^3 \wedge e^4},\, \\ & \qquad {a_{13}a_{24}\ne 0}\\
{\rr_{4,-1, \beta}:} & \omega={a_{12}e^1 \wedge e^2 + a_{14}e^1 \wedge e^3+a_{24}e^2 \wedge e^4 + a_{34}e^3 \wedge e^4},\, \\ & \qquad {a_{14}a_{23}\ne 0},\,  {\beta \ne -1,0,1} \\
{\rr_{4,-1, -1}:} & \omega={a_{12}e^1 \wedge e^2 + a_{13}e^1 \wedge e^3+a_{14}e^1 \wedge e^4+a_{24}e^2 \wedge e^4+a_{34}e^3 \wedge e^4 },\,\\ &
\qquad  {a_{12}a_{34} - a_{13}a_{24}\ne 0} \\
{\rr_{4,\alpha, -\alpha}:} & \omega={a_{14}e^1 \wedge e^4 + a_{23}e^2 \wedge e^3+a_{24}e^2 \wedge e^4 +a_{34}e^3 \wedge e^4},\,\\
& \qquad  {a_{14}a_{23}\ne 0},\, {\alpha \ne -1,0}\\
{\rr'_{4,0, \delta}:} & \omega={a_{14}e^1 \wedge e^4 + a_{23}e^2 \wedge e^3+a_{24}e^2 \wedge e^4 +a_{34} e^3 \wedge e^4},\, \\ & \qquad {a_{14}a_{23}\ne 0},\, \delta \ne 0\\
{\dd_{4,1}:} & { \omega=a_{12-34}(e^1 \wedge e^2 - e^3 \wedge e^4) +a_{14}e^1 \wedge e^4 + a_{24}e^2 \wedge e^4},\,\\ & \qquad  {a_{12-34}\ne 0} \\
{\dd_{4,2}:} & {\omega=a_{12-34} (e^1\wedge e^2-e^3 \wedge e^4)+a_{14}e^1 \wedge e^4 +a_{23} e^2 \wedge e^3+ a_{24} e^2 \wedge e^4},\,\\ & \qquad  {-a_{12-34}^2+a_{14}a_{23}\ne 0} \\
{\dd_{4,\lambda}:} & {\omega=a_{12-34} (e^1\wedge e^2-e^3 \wedge e^4)+a_{14}e^1 \wedge e^4 +a_{24} e^2 \wedge e^4},\,\\
& \qquad  {a_{12-34}\ne 0},\, \lambda  \ne 1,2\\
{\dd'_{4,\delta}:} & {\omega=a_{12-\delta34} (e^1\wedge e^2-\delta e^3 \wedge e^4)+a_{14}e^1 \wedge e^4 +a_{24} e^2 \wedge e^4}\\
& \qquad  {a_{-12+\delta 34}\ne 0},\, \delta  \ne 0\\
{\hh_4:} & {\omega=a_{12-34} (e^1\wedge e^2-e^3 \wedge e^4)+a_{14}e^1 \wedge e^4 +a_{24} e^2 \wedge e^4} \\
& \qquad  {a_{12-34}\ne 0} \\
\end{array}
$$
\end{prop}
\begin{proof} The proof follows by working on each Lie algebra and searching on it for closed two-forms. Thus let $\omega$ be a two-form, that is $\omega = \sum_{i\ge 1, j>i} e^i \wedge e^j$. Let $d$ be the antiderivation operator, using that $d(e^i \wedge e^j) = de^i \wedge e^j - e^i\wedge de^j$ and the Lie bracket relations of Proposition \ref{clases} one determines closed two forms. The next step is to  find the rank of $\omega$, that is to compute $\omega^2$. If $\omega$  has maximal rank, then $\ggo$ will be endowed with a symplectic structure. As all cases should be handled in a similar way, we will give in detail the computations on  $\rr'_2$,  the Lie algebra which corresponds to $\aff (\CC)$.

  The Lie bracket relations on $\ggo$ implies that $d e^1 =0=de^2$ and  $-de^3 = e^1\wedge e^3-e^2\wedge e^4$, and $-de^4 = e^1 \wedge e^4 + e^2 \wedge e^3$. At the next level we have:
$$ d(e^1 \wedge e^3) = -e^1\wedge e^2 \wedge e^4,\, d(e^1 \wedge e^4) = e^1\wedge e^2 \wedge e^3,\, d(e^2 \wedge e^3) = -e^1\wedge e^2 \wedge e^3$$
$$d(e^2 \wedge e^4) =- e^1\wedge e^2 \wedge e^4,\qquad d(e^3 \wedge e^4) = -2e^1\wedge e^3 \wedge e^4$$
Thus any 2-form $\omega$ that is closed has the form $\omega = a_{12}e^1 \wedge e^2+a_{13-24}(e^1\wedge e^3- e^2 \wedge e^4) +a_{14+23} (e^1\wedge e^4+e^2 \wedge e^3)$. Now $\omega$ is symplectic if it satisfies the conditions mentioned in the Table (\ref{simplecticas}) for $\rr'_2$, concluding the proof.
\end{proof}

\medspace

Note that among the symplectic four dimensional Lie algebras the unimodular ones are $\rr\rr_{3,-1}$, $\rr\rr_{3,0}'$ and $\nn_4$.

\medspace

Recall that an element $\omega \in \Lambda^p(\ggo^{\ast})$ is called exact if $\omega = d \eta$ for some $\eta \in \Lambda^{p-1}(\ggo^{\ast})$. Thus the computations of the previous propositions give  also the exact symplectic Lie algebras, which were obtained by Campoamor in \cite{Ca}.

\begin{cor} A four dimensional solvable Lie algebra admits an exact symplectic structure if and only if $\ggo$ is one of the following attached with the respective symplectic structure

\begin{center}
\small{
\begin{tabular}{|c|c|c|}\hline
{{\rm Case}} & {{$\omega$}} & {\rm Condition}\\ \hline
{$\rr_2\rr_2$} & {$a_{12}e^1 \wedge e^2 + a_{34}e^3 \wedge e^4$} & {$a_{12}a_{34}\ne 0$} \\\hline
{$\rr'_2$} & {$ a_{13-24}(e^1 \wedge e^3 - e^2 \wedge e^4) +a_{14+23}(e^1 \wedge e^4+e^2 \wedge e^3)  $} & {$a_{14+23}^2 + a_{13-24}^2\ne 0$} \\ \hline
{$\dd_{4,1}$} & { $a_{12-34}(e^1 \wedge e^2 - e^3 \wedge e^4) +a_{14}e^1 \wedge e^4$} & {$a_{12-34}\ne 0$} \\ \hline
{$\dd_{4,\lambda}\lambda \ne 1$} & { $a_{12-34}(e^1 \wedge e^2 - e^3 \wedge e^4) +a_{14}e^1 \wedge e^4 + a_{24}e^2 \wedge e^4$} & {$a_{12-34}\ne 0$} \\ \hline
{$\dd'_{4,\delta}\delta \ne 0$} & {$a_{-12+\delta34} (-e^1\wedge e^2+\delta e^3 \wedge e^4)+a_{14}e^1 \wedge e^4 +a_{24}e^2 \wedge e^4$} & {$a_{-12+\delta 34}\ne 0$} \\\hline
{$\hh_4$} & {$a_{12-34} (e^1\wedge e^2-e^3 \wedge e^4)+a_{14}e^1 \wedge e^4 + a_{24} e^2 \wedge e^4$} & {$a_{12-34}\ne 0$} \\
\hline
\end{tabular} }
\end{center}
\end{cor}

Recall that two symplectic Lie algebras $(\ggo_1,\omega_1)$ and $(\ggo_2,\omega_2)$ are said to be symplectomorphically equivalent if there exists an isomorphism of Lie algebras $\varphi:\ggo_1 \to \ggo_2$, which preserves the symplectic forms, that is $\varphi^{\ast}\omega_2 = \omega_1$.

\begin{prop} Let $\ggo$ be a symplectic real Lie algebra of dimension four. Then $\ggo$ is symplectomorphically equivalent to one of the following Lie algebras equipped with a symplectic form as follows:

$$
\begin{array}{ll}
{\rr\hh_{3}:} & {\omega=e^1 \wedge e^4 + e^2 \wedge e^3} \, \\
{\rr \rr_{3,0}:}& \omega=e^1 \wedge e^2 + e^3 \wedge e^4\,\\
{\rr \rr_{3,-1}:}& \omega= e^1 \wedge e^4 + e^2 \wedge e^3\, \\
{\rr\rr'_{3,0}:} & \omega=e^1 \wedge e^4 + e^2 \wedge e^3 \,\\
{\rr_2\rr_2:} & \omega_{\lambda}=e^1 \wedge e^2 + \lambda e^1 \wedge e^3 +e^3 \wedge e^4,\quad \lambda \ge 0\\
{\rr'_2:} & \omega={e^1 \wedge e^4+e^2 \wedge e^3 } \,\\
{\nn_{4}:} & \omega={e^1 \wedge e^2 +e^3 \wedge e^4}\, \\
{\rr_{4,0}:} & \omega_+={e^1 \wedge e^4 + e^2 \wedge e^3,\quad \omega_-=e^1 \wedge e^4 - e^2 \wedge e^3 }\, \\
{\rr_{4,-1}:} & \omega=e^1 \wedge e^3 +e^2 \wedge e^4 \\
{\rr_{4,-1, \beta}:} & \omega={e^1 \wedge e^2 + e^3 \wedge e^4},\quad    {-1 \le  \beta < 1} \\
{\rr_{4,\alpha, -\alpha}:} & \omega={e^1 \wedge e^4 + e^2 \wedge e^3, \quad -1 < \alpha <0}\\
{\rr'_{4,0, \delta}:} & \omega_+={e^1 \wedge e^4 + e^2 \wedge e^3,\quad \omega_-=e^1 \wedge e^4 - e^2 \wedge e^3 },\quad  \delta > 0\\
{\dd_{4,1}:} & { \omega_1=e^1 \wedge e^2 - e^3 \wedge e^4,\quad  \omega_2=e^1 \wedge e^2 - e^3 \wedge e^4 +e^2 \wedge e^4}\,\\
{\dd_{4,2}:} & {\omega_1=e^1\wedge e^2-e^3 \wedge e^4,\, \omega_2=e^1 \wedge e^4 + e^2 \wedge e^3,\, \omega_3= e^1\wedge e^4-e^2 \wedge e^3}\\
{\dd_{4,\lambda}:} & {\omega=e^1\wedge e^2-e^3 \wedge e^4},\quad  \lambda\ge \frac12, \quad \lambda  \ne 1,2\\
{\dd'_{4,\delta}:} & {\omega_+=e^1\wedge e^2-\delta e^3 \wedge e^4, \quad  \omega_-=-e^1\wedge e^2+\delta e^3 \wedge e^4},\quad  \delta >0\\
{\hh_4:} & {\omega_+=e^1\wedge e^2-e^3 \wedge e^4,\quad  \omega_-=-e^1\wedge e^2+e^3 \wedge e^4} \\
\end{array}
$$
\end{prop}
\begin{proof} The proof make use of the authomorphisms of symplectic Lie algebras computed in Appendix I. As an example, for the case $\rr_2'$ the authomorphism given by $\sigma e_1= e_1$, $\sigma e_2 = e_ 2 + \alpha e_4$, $\sigma e_3 =  \gamma e_3 - \beta e_4$, $\sigma e_4 = \beta e_3 + \gamma e_4$ does $\sigma^{\ast}(e^1 \wedge e^4+ e^2 \wedge e^3) = {\alpha e^1 \wedge e^2 + \beta (e^1 \wedge e^3 - e^2 \wedge e^4) +\gamma (e^1 \wedge e^4+e^2 \wedge e^3) }$ where $\beta^2 + \gamma^2 \neq 0$. Similar computations on each symplectic Lie algebra complete the proof.
\end{proof}
\

 In the following   we  relate the existence problem of symplectic forms in the four dimensional case  with  the action of the adjoint representation of the corresponding Lie algebra. That is, we read the list of symplectic Lie algebras of Proposition \ref{simplecticas} from the point of view of the adjoint action of the Lie algebra. That gives rise to obstructions to the existence of symplectic structures in  higher dimensions  (see  section \ref{last}).

\begin{cor}\label{sympuni} Let $\ggo$ be a unimodular  four dimensional solvable non abelian Lie algebra. Then $\ggo$ is symplectic if and only if $\ggo$ is isomorphic to $\nn_4$ or $\ggo$ is isomorphic to a direct product of $\RR$ and a three dimensional unimodular solvable Lie algebra.
\end{cor}

\begin{remark}  In \cite{LM} the authors proved that unimodular symplectic Lie algebras must be solvable.
\end{remark}

\begin{cor} \label{sympnouni} Let $\ggo$ be a non unimodular  four dimensional solvable  Lie algebra. If $\ggo$ is symplectic then either:

$\ggo' \simeq \RR$ or

$\ggo' \simeq \RR^2$ and if $\ggo$ is not isomorphic to $\rr_{4,0}$, then $\zz(\ggo)$ is trivial, that is $\ggo$ is isomorphic either to $\rr_2\rr_2$ or $\rr_2'$, or

$\ggo' \simeq \hh_3$ or

$\ggo' \simeq \RR^3$ and the adjoint action of an element $e_0 \notin \ggo'$ is equivalent to one of the following ones:
$$
\begin{array}{cccc}
\left(
\begin{matrix}
1 & 0 & 0 \\
0 & -1 & 1 \\
0 & 0 & -1
\end{matrix}
\right) &
\left(
\begin{matrix}
1 & 0 & 0 \\
0 & -1 & 0 \\
0 & 0 & \beta
\end{matrix}
\right) & \left(
\begin{matrix}
1 & 0 & 0 \\
0 & \alpha & 0 \\
0 & 0 & -\alpha
\end{matrix}
\right) &
\left(
\begin{matrix}
1 & 0 & 0 \\
0 & 0 & \delta \\
0 & -\delta & 0
\end{matrix}
\right) \\
& -1\le \beta < 0 \,& \, -1<\alpha < 0 \,& \\
\end{array}
$$
\end{cor}

\section{Models for symplectic Lie algebras}

Models for symplectic Lie algebras were given by different authors. In this section we describe symplectic four dimensional Lie algebras in terms of two basic constructions: either as solutions of the cotangent extension problem or as symplectic double extensions. Essentially the difference between the two cases is the existence of lagrangian ideals in the first one and  of  isotropic non lagrangian ideals in the second one.

\vspace{.3cm}

Let  $(\ggo, \Omega)$  be a Lie algebra endowed with a non-degenerate skew-symmetric bilinear form. If  $W\subset \ggo$ is a subspace of $\ggo$ then the orthogonal subspace $W^{\perp}$ is
$$W^{\perp} = \{ x \in \ggo \, /\, \Omega(x,y) =0 \, \text{ for all } y \in W\}$$
In particular it always holds that $\dim \ggo = \dim W + \dim W^{\perp}$. The subspace $W$ is
 called {\it  isotropic} if $\Omega (W,W) = 0$ (that is $W\subset W^{\perp}$) and is called {\it lagrangian} if $W^{\perp}=W$.

It is easy to see that a subspace $W$ is lagrangian if and only if $W$ is isotropic and $\dim \ggo = 2 \dim W$. Moreover, since $\Omega$ is closed, an isotropic ideal $W$ must be abelian and $W^{\perp}$ is a subalgebra.

\begin{lem}\label{isotropic} Let $\ggo$ be a  symplectic four dimensional  Lie algebra, then $\ggo$ always admits an isotropic ideal $\jj$. Moreover except for the Lie algebras $\rr\rr_{3,0}',\rr_{4,0,\delta}', \dd_{4,\lambda}'$ all other Lie algebras admit  lagrangian ideals.
\end{lem}
\begin{proof} For each symplectic four dimensional Lie algebra we will exhibit an isotropic ideal with respect to every symplectic form - see Proposition \ref{simplecticas}.
$$
\begin{array}{llllll}
\rr\hh_3: & \la e_3, e_4\ra,\quad & \rr\rr_{3,-1}: & \la e_2, e_4\ra \quad &  \rr\rr_{3,0}: & \la e_2, e_3\ra,\quad \\
 \rr\rr'_{3,0}: & \la e_4\ra, & \rr_2\rr_2: & \quad \ggo',\quad  & \rr_2': & \quad \ggo'\quad \\
\nn_4: & \quad \ggo',\quad  & \rr_{4,0}: & \quad \ggo' &\rr_{4,-1}: & \la e_1, e_2\ra, \\
 \rr_{4,-1,-1}: & \la e_2, e_3\ra  &\rr_{4,-1,\beta}: & \la e_1, e_2\ra, & \rr_{4,\alpha,-\alpha}: & \la e_1, e_3\ra, \\
 \rr'_{4,0,\delta}: & \la e_1\ra, & \dd_{4,\lambda}: & \la e_1, e_3\ra, & \dd'_{4,\delta}: & \la e_3\ra, \\
\hh_{4}: & \la e_1, e_3\ra
\end{array}
$$
Moreover $\jj \subset \ggo'+\zz(\ggo)$ and is abelian. 
\end{proof}

\subsection{Cotangent extension problem}
Let $\hh$ be a Lie algebra and let $\hh^{\ast}$ be
the dual vector space of $\hh$, consider $\hh^{\ast} \oplus \hh$  as a vector space and let $\omega_0$ be the skew symmetric two form  defined as:
\begin{equation}\label{omega}
\omega_0((\varphi, x), (\varphi', x')) = \varphi(x') -\varphi'(x)\end{equation}
 The {\it cotangent extension problem} consists in finding a Lie algebra structure on $\hh^{\ast} \oplus \hh$ which  satisfy the following conditions:

(c1) $0 \longrightarrow \hh^{\ast} \longrightarrow \hh^{\ast} \oplus \hh \longrightarrow \hh \longrightarrow 0$ is an exact sequence of Lie algebras, $\hh^{\ast}$ endowed with its abelian Lie algebra structure;

(c2) the left invariant 2-form spanned by $\omega_0$ is closed.

A symplectic Lie algebra $(\ggo, \omega)$ is said to be a solution of the cotangent extension problem if $\ggo$ is symplectomorphically equivalent to a Lie algebra of the form $(\hh^{\ast} \oplus \hh, \omega_0)$ satisfying conditions (c1) and (c2). Thus $\ggo$ is an extension of the Lie algebra $\hh$.

In what follows we describe conditions to get solutions of the cotangent extension problem.

\

Let ($\hh, [\, , \,]_{\hh})$ be a Lie algebra and $\rho: \hh \to \End (\hh^{\ast})$ be a representation of $\hh$ on the dual space of left invariant 1-forms of $\hh$, denoted $\hh^{\ast}$. Thus $\hh^{\ast}$ inherits a structure of $\hh$-module, denoted $x . \varphi = \rho(x) \varphi$.

\medspace

 Let $\ggo$ be the direct sum of the vector spaces $\hh^{\ast} \oplus \hh$ and define a skew-symmetric map on $\ggo$,
 $[\, , \,]: \ggo \times \ggo \to \ggo$ by:
$$
\begin{array}{llll}
 \left[\varphi, \eta \right] & = & 0 & \varphi, \, \eta \in \hh^{\ast} \\
\left[ \varphi, x  \right] & = & - x . \varphi & x \in \hh,\,\varphi \in \hh^{\ast} \\
 \left[ x,y \right] & = & [x,y]_{\hh} + \alpha(x,y) & x, \, y \in \hh
\end{array}
$$
where $\alpha$ is a 2-cochain. Then $[ , ]$ defines a structure of Lie algebra on $\ggo$ if and only if $\alpha$ is a 2-cocycle,  $\alpha \in Z^2(\hh, \hh^{\ast})$, that is
\begin{equation}
\label{lie}
\alpha([x_1, x_2]_{\hh}, x_3) + \alpha([x_2, x_3]_{\hh}, x_1)+ \alpha([x_3, x_1]_{\hh}, x_2) = \end{equation}
\[ \qquad  x_3.\alpha(x_1,x_2) +  x_1.\alpha(x_2,x_3) +  x_2.\alpha(x_3,x_1)
\]

In this situation $\ggo$ is an extension of $\hh$ and  we have the following exact sequence of Lie algebras
\[ 0 \longrightarrow \hh^{\ast} \longrightarrow \ggo \longrightarrow \hh \longrightarrow 0
\]

\medspace

Denoting by $<,>$ the evaluation operation on $\hh^{\ast} \oplus \hh$ then the two-form $\omega_0$ on $\ggo$ defined by
$$\omega_0((\varphi_1,x_1),(\varphi_2,x_2)) = <\varphi_1,x_2> -<\varphi_2,x_1>$$
is closed if and only if
\begin{equation}
\label{sum}
<\alpha(x_1,x_2),x_3>+ <\alpha(x_2,x_3),x_1> + <\alpha(x_3,x_1),x_2> = 0
\end{equation} for all $x_1, x_2, x_3 \in \hh$ and
\begin{equation}
\label{coborde}
<x.\varphi, y> -< y. \varphi,x> = <\varphi,[x,y]_{\hh}>
\end{equation}

The  condition (\ref{sum}), known as the ``Bianchi identity'', is equivalent to say that $\alpha$ belongs to the kernel of the canonical map $(\Lambda^2 \hh \otimes \hh \to \Lambda^3 \hh)$, which is $\mathbb{S}_{(2,1)}(\hh)$ the Weyl space corresponding to the partition 3=2+1  (see \cite{FH}).

Then the  resulting Lie algebra $\ggo$ (attached to the  triple $(\hh, \rho, [\alpha])$) satisfying (\ref{lie}), (\ref{sum}) and (\ref{coborde}) is a  {\it solution of the cotangent extension problem}.

\

\begin{remark} \label{cano} Note that the coadjoint representation satisfies (\ref{coborde}) if and only if $\hh$ is abelian. Thus $\omega_0$ spans a symplectic structure on the cotangent bundle $T^{\ast}H$ of a Lie group $H$, endowed with its canonical Lie group structure, if and only if $H$ is abelian.
\end{remark}

\begin{remark}
If one defines a connection $\nabla$ on $\hh$ by
$$<\varphi, \nabla_xy> = -<x.\varphi,y> \qquad (\mbox{so } \varphi \circ \nabla_x = - x. \varphi) \quad x,y \in \hh, \,\varphi \in \hh^{\ast}$$
then (\ref{coborde}) is equivalent to the fact that the connection is torsion free (see  \cite{By1}).
\end{remark}

\begin{remark} Observe that $\hh$ is a subalgebra of $\ggo$ if and only if $\alpha=0$ and in this case $\ggo$ is the semidirect product of $\hh$ and $\hh^{\ast}$. Moreover  $\omega_0$ is closed if and only if (\ref{coborde}) is verified. This case was studied by Boyom in \cite{By1}.
\end{remark}

\begin{remark} If the representation $\rho$ is trivial, then (\ref{lie}) becomes
$$
\alpha([x_1, x_2]_{\hh}, x_3) + \alpha([x_2, x_3]_{\hh}, x_1)+ \alpha([x_3, x_1]_{\hh}, x_2) = 0
$$
 and $\omega_0$ is closed if and only if (\ref{sum}) holds and (\ref{coborde}) becomes
$$
 <\varphi , [x,y]_{\hh}> = 0
$$
So $\hh$ must be abelian (direct sum of vector spaces) and  $\ggo$ is a two-step nilpotent Lie algebra with Lie bracket $[,]$ defined by
$$[x,y] = \alpha(x,y) \quad x,\, y \in \hh$$
and so $\hh' = \Im \alpha$.  The 2-form $\omega_0$ is closed if and only if  (\ref{sum}) holds.
\end{remark}

The previous explanation proofs the first assertion of the following theorem. The second assertiont relates Lie algebras having a lagrangian ideal with solutions of the cotangent extension problem. We include the proof, which is  useful for our purposes. However similar results are known in a more general context (see (\ref{We2})).

\begin{thm} \label{cota} i) Let $\ggo$ be a 2n-dimensional Lie algebra with an abelian ideal $\hh^{\ast}$ of dimension n (so we have (c1)). Then $\ggo$ is a solution of the cotangent extension problem, if and only if conditions  (\ref{sum}) and (\ref{coborde}) are satisfied.

ii)  Let $(\ggo, \omega)$ be a symplectic Lie algebra with a lagrangian ideal. Then $\ggo$ is a solution of the cotangent extension problem.
\end{thm}
\begin{proof} ii)  Let $\jj$ be a lagrangian ideal, then one has the following exact sequences of Lie algebras:
$$0 \longrightarrow \jj \longrightarrow \ggo \longrightarrow \ggo/\jj \longrightarrow 0$$
Let $\hh:=\ggo/\jj$ be the quotient Lie algebra. Since $\jj$ is abelian it can be identified  with $\hh^{\ast}$, the last one endowed with the abelian Lie algebra structure.  Let $\ggo= \jj \oplus \vv$ be a splitting into lagrangian subspaces (this always exists, see \cite{We} Lect. 2), then the map $\beta: \jj \to (\hh)^{\ast}$ given by $\beta(x)(y+\jj) = \omega(x,y)$ for $x\in  \jj$, $y+\jj \in h=\ggo/\jj$, induces an isomorphism $\tilde{\beta}:\jj \to \vv^{\ast}$. Giving $\hh^{\ast}\oplus \hh$ the Lie algebra structure via the isomorphism $1 \oplus \tilde{\beta}: \jj \oplus \vv \to \hh^{\ast} \oplus \hh$, one may check that $1 \oplus \tilde{\beta}$  is a symplectomorphism from $(\ggo, \omega)$ to $(\hh^{\ast} \oplus \hh, \omega_0)$, completing the proof of the second assertion of the theorem.
\end{proof}

\begin{remark} \label{We2} It is known that every lagrangian foliation is locally symplectomorphic to the foliation of $\RR^{2n}$ by the manifolds $x_i=$constant, and that the leaves of a lagrangian foliation carry a natural flat torsion free affine connection (see \cite{We2}).
\end{remark}

\begin{defn} Two solutions of the cotangent extension problem $\ggo_1$ and $\ggo_2$ resulting as extensions of $\hh_1$ and $\hh_2$ respectively are said to be equivalent if there exists a isomorphism $\psi:\ggo_1 \to \ggo_2$ such that  $\psi \hh_1 = \hh_2$
\end{defn}

\begin{thm} The following tables show all solutions $\ggo$ up to equivalence of the cotangent extension problem for $\hh$ of dimension two.\end{thm}
The elements of  $H^2_{\rho}(\hh, \RR^2)$ are denoted by $\alpha$.  The basis of $\hh$ is $\{x,y\}$ and for $\hh^{\ast}$ we choose the corresponding dual basis. The symbol (*) indicates that the  Lie algebra obtained as an extension of $\hh$, via $\rho$ and $\alpha$, is a solution of the cotangent extension problem.

\begin{center}
\small{
\begin{tabular}{|c|c|c|}\hline
{{\rm Representation $\rho$}} & {{$H^2_{\rho}(\RR^2, \RR^2)$}} & {$\ggo$} \\ \hline
{$\rho \equiv 0$} & {2} & {$\left\{ \begin{array}{ll}
\RR^4 & \alpha = 0 (*)\\
\rr \hh_3 & \alpha \neq 0 (*)
\end{array} \right. $} \\ \hline
{$\begin{array}{ll}
\rho(x) = 0 & \rho(y) = \left( \begin{matrix} 0 & 0 \\ 0 & 1 \end{matrix} \right)\end{array}$} &
1  & {\begin{tabular}{c}
{$\left\{ \begin{array}{ll}
\rr \rr_{3,0} & \alpha = 0 (*) \\
\rr_{4,0} & \alpha \neq 0 (*) \end{array}\right.$}\end{tabular}}\\ \hline
{$\begin{array}{ll}
\rho(x) = 0 & \rho(y) = \left( \begin{matrix} \lambda & 0 \\ 0 & 1 \end{matrix} \right)\\
& \qquad \lambda \neq 0\end{array}$} &
0 & {$
\rr \rr_{3,\lambda} $}
\\ \hline
{$\begin{array}{ll}
\rho(x) = 0 & \rho(y) = \left( \begin{matrix} 0 & 0 \\ 1 & 0 \end{matrix} \right)\end{array}$} & {1}  &
{\begin{tabular}{c}
{$\left\{ \begin{array}{ll}
\rr \hh_{3} & \alpha = 0 (*)\\
\nn_{4} &  \alpha \neq 0 (*) \end{array} \right.$}\end{tabular}}
\\ \hline
{$\begin{array}{ll}
\rho(x) = 0 & \rho(y) = \left( \begin{matrix} \lambda & 0 \\ 1 & \lambda \end{matrix} \right)\\
 & \qquad \lambda\neq 0 \end{array}$} & {
0 }  &
{$\rr \rr_{3} $}
\\ \hline
{$\begin{array}{ll}
\rho(x) = 0 & \rho(y) = \left( \begin{matrix} \gamma & 1 \\ -1 & \gamma \end{matrix} \right)\end{array}$} & {0 } & {$\rr \rr_{3,\gamma}' $}\\ \hline
{$\begin{array}{ll}
\rho(x) = \left( \begin{matrix} 1 & 0 \\ 0 & 0 \end{matrix} \right) & \rho(y) = \left( \begin{matrix} 0 & 0 \\ 0 & 1 \end{matrix} \right)\end{array}$} & {0} & {$\rr_2 \rr_2 \quad(*) $}\\ \hline
{$\begin{array}{ll}
\rho(x) = \left( \begin{matrix} 1 & 0 \\ 0 & 1 \end{matrix} \right) & \rho(y) = \left( \begin{matrix} 0 & 1 \\ 0 & 0 \end{matrix} \right)\end{array}$} & {0} & {$\dd_{4,1} \quad(*) $}\\ \hline
{$\begin{array}{ll}
\rho(x) = \left( \begin{matrix} 1 & 0 \\ 0 & 1 \end{matrix} \right) & \rho(y) = \left( \begin{matrix} 0 & -1 \\ 1 & 0 \end{matrix} \right)\end{array}$} & {0} & {$\rr_2' \quad(*) $}\\ \hline\end{tabular}}
\end{center}
\begin{center}Table \label{rr2} \ref{rr2} for $\hh= \RR^2$\end{center}


\begin{center}
\small{
\begin{tabular}{|c|c|c|}\hline
{{\rm Representation $\rho$}} & {{$H^2_{\rho}(\aff(\RR), \RR^2)$}} & {$\ggo$} \\ \hline
{$\rho \equiv 0$} & {0} & {$\rr \rr_{3,0}(*)
 $} \\ \hline
{$\begin{array}{ll}
\rho(y) = 0 & \rho(x) = \left( \begin{matrix} 0 & 0 \\ 0 & \lambda \end{matrix} \right)\\
& \lambda \in [-1,1] \end{array}$} & {0} & {$
\rr \rr_{3,\lambda} (**)$}\\ \hline
{$\begin{array}{ll}
\rho(y) = 0 & \rho(x) = \left( \begin{matrix} \mu & 0 \\ 0 & \lambda \end{matrix} \right)\\
 -1 < \mu \leq \lambda, \mu \lambda \neq 0 & \text{ or }
-1 = \mu \leq \lambda \leq 0
\end{array}$} & {0} & {$
\rr_{4,\mu,\lambda} (***) $}\\ \hline
{$\begin{array}{ll}
\rho(y) = 0 & \rho(x) = \left( \begin{matrix} \mu & 1 \\ 0 & \mu \end{matrix} \right)\end{array}$} & {0}  & $\rr_{4,\mu}$ \\ \hline
{$\begin{array}{ll}
\rho(y) = 0 & \rho(x) = \left( \begin{matrix} \gamma & \delta \\ -\delta & \gamma \end{matrix} \right)\\ & \gamma \in \RR,\, \delta > 0 \end{array}$} & {0 } & {$\rr_{4,\gamma, \delta}' $}\\  \hline
{\begin{tabular}{c}
{$
\rho(y) = \left( \begin{matrix} 0 & 0 \\ 1 & 0 \end{matrix} \right)$} \\
{$\rho(x) = \left( \begin{matrix} 1 & 0 \\ 0 & 2 \end{matrix} \right)$}\end{tabular}} &
 1  & {$\left\{\begin{array}{ll}
\dd_{4,1/2} & \alpha =  0 (*)\\
\hh_{4} &  \alpha \neq 0 (*)
\end{array} \right.$}\\\hline
{\begin{tabular}{c}
{$\rho(y) = \left( \begin{matrix} 0 & 0 \\ 1 & 0 \end{matrix} \right)$} \\
{$\rho(x) = \left( \begin{matrix} -1 & 0 \\ 0 & 0 \end{matrix} \right)$}  \end{tabular}} &
{0 } &
{$\dd_{4} $}\\
 \hline
{\begin{tabular}{c}
{$\rho(y) = \left( \begin{matrix} 0 & 0 \\ 1 & 0 \end{matrix} \right)$} \\
{$\rho(x) = \left( \begin{matrix} a - 1/2 & 0 \\ 0 & \ a +1/2 \end{matrix} \right)$} \\  \qquad {$ a\neq 3/2,-1/2 $} \end{tabular}} &
{0 } &
{$
\dd_{4,\frac{1}{a+1/2}}  (*)
$}\\
 \hline

\end{tabular}}
\end{center}
\begin{center}Table \label{aff} \ref{aff} for $\hh= \aff(\RR)$ \end{center}

\medspace

(**) solution for $\lambda = -1$.

(***) solution for $(\mu, \lambda)$ of the form  (-1,-1) $(-1, \lambda),(\mu,-\mu)$ with restrictions of Proposition \ref{simplecticas}.
\begin{remark} In \cite{A-B-D-O} it is was considered a special case of the sequence of Lie algebras (c1), that is when  $\hh$ is a two-dimensional Lie algebra and the exact sequence (c1) splits, so that $\ggo$ is a semidirect product of $\hh$ and $\RR^2$. \end{remark}
\begin{proof}
To construct the tables, we make use of the information given in  \cite{A-B-D-O} to get all semidirect  extensions of $\hh$. The idea of this proof is to study the image of  $\rho$ in $\glo(2,\RR)$. Thus if $\hh \simeq \RR^2$, then the image of $\rho$ can be described in terms of the following  subalgebras of $\glo(2,\RR)$:
$$
\left\{
\left( \begin{matrix} a & 0 \\
0 & b \end{matrix} \right): a, b \in \RR \right\}, \quad
\left\{
\left( \begin{matrix} a & 1 \\
0 & a \end{matrix} \right): a \in \RR \right\}, \quad
\left\{
\left( \begin{matrix} a & b \\
b & -a \end{matrix} \right): a, b \in \RR \right\}
$$
If the image of $\rho$ is one dimensional, then $\ggo$ is an extension of a three dimensional Lie algebra. If the image of $\rho$ is two-dimensional, then $\ggo$ is either isomorphic to $\rr_2 \rr_2$ or to $\rr_2'$.

If $\hh$ is isomorphic to $\aff(\RR)$, then we may assume that $\hh = \la x, y\ra$, with $[x,y]=y$. Thus if $\rho$ is trivial then $\ggo\simeq \rr \rr_{3,0}$. If the image of $\rho$ is one dimensional then $\rho(y)=0$ and $\rho(x)$ acts on $\hh^{\ast}$ as follows:
$$
\left\{
\left( \begin{matrix} \mu & 0 \\
0 & \lambda \end{matrix} \right): \mu, \lambda \in \RR, \lambda \neq 0 \right\}, \quad
\left\{
\left( \begin{matrix} \mu & 0 \\
1 & \mu \end{matrix} \right): \mu \in \RR \right\}, \quad
\left\{
\left( \begin{matrix} \gamma & \delta \\
\delta & -\gamma \end{matrix} \right): \gamma, \delta \in \RR, \delta \neq 0 \right\}
$$
In these cases we get the first part of the second table. If the image of $\rho$ is two dimensional, then one may assume that
$$\rho(y) = \left(
\begin{matrix} 0 & 0\\
1 & 0 \end{matrix} \right)
$$
and $\rho(x)$ takes the following form
$$ \rho(x) = \left(
\begin{matrix} \alpha +1/2 & 0\\
0 & \alpha-1/2 \end{matrix} \right)
$$
and so we get all semidirect extensions. To complete the proof we need to compute the second cohomology group and to determine the resulting Lie algebra. This computations give the results on the second column. In this way we have completed the work of \cite{A-B-D-O} by obtaining all extensions in the two dimensional case, up to equivalence. Comparing with notations of Proposition (\ref{clases}) we conclude the proof.
\end{proof}

\begin{thm}\label{model} Let  $\ggo$ be four dimensional Lie algebra; if $\ggo$  is a solution of the cotangent extension problem then $\ggo$ is either sympletic completely solvable or isomorphic to $\aff(\CC)$ \end{thm}

\begin{proof}  A first proof follows by reading the results of the previous tables.

A second proof is obtained as follows: by  Lemma \ref{isotropic} there always exists a abelian lagrangian ideal $\jj$ on any symplectic  completely solvable four dimensional Lie algebra. The proof will be completed be applying Theorem (\ref{cota}).
\end{proof}

\begin{remark} The symplectic four dimensional  Lie algebras $\rr\rr_{3,0}',\rr_{4,0,\delta}'$ do not admit lagrangian ideals  however they are semidirect products of two symplectic 2-dimensional Lie algebras. In fact in both cases  take the symplectic subalgebras $\hh=<e_1,e_4>$ and the abelian one $<e_2,e_3>$ endowed with the induced symplectic structure of $\ggo$. The Lie algebra $\hh$ is abelian for $\rr\rr_{3,0}'$ and is isomorphic to $\aff(\RR)$ for $\rr_{4,0,\delta}'$. But $\dd_{4,\delta}'$ cannot be written as a semidirect product of  2-dimensional Lie algebras. It was proved in \cite{A-B-D-O} that this Lie algebra does not admit a decomposition as a direct sum of two 2-dimensional subalgebras (sum as vector spaces).
\end{remark}

\subsection{Symplectic double extensions \cite{MR}\cite{DM1}}


In this section we modelize some four dimensional symplectic Lie algebras as symplectic double extensions in a classical or in  a generalized sense. To this end we  recall the main ideas of these constructions of symplectic Lie algebras, and we remit to the papers of Medina A. and Revoy P. \cite{MR} or Dardi\' e J. and Medina A. \cite{DM1} for the details.

Let (B, $\omega'$) be a symplectic Lie algebra, let $\delta$ be a derivation of B and let $z \in$ B. Let $I = \RR e \oplus B$ be the central extension of B by $\RR e$ defined by
$$ [ a, b]_I = [a,b]_B + \omega'(\delta a , b) e\qquad a,b \in B$$
where $[\,, \, ]_B$ denotes the Lie bracket on $B$. Let $A$ be the semidirect product of $I$ by $\RR d$ given by
$$[d, e] = 0 \qquad [d,a] = -\omega'(z,a)e -\delta(a) \qquad a \in B$$
Extending the symplectic structure of B to $A$ via $\omega$, defined as $\omega(e,d) = 1$, and $\omega(B, e) = 0 = \omega(B,d)$, then $A$ is said to be a symplectic double extension of (B, $\omega')$ by $\RR$.

\begin{thm} \cite{MR} Let (A, $\omega$) be a symplectic  2n-dimensional Lie algebra with non trivial center. Then A is a symplectic double extension by $\RR$ of a 2n-2 dimensional symplectic Lie algebra (B,$\omega')$.
\end{thm}

This result motivates the following definition in \cite{DM1}.

\begin{defn} A symplectic Lie algebra $(\ggo,\omega)$ is called a symplectic double extension of a symplectic Lie algebra $(W,\omega')$ if there exists a central one-dimensional subalgebra $\jj\subset \ggo$ such that the reduced symplectic Lie algebra  $\jj^{\perp}/\jj$ is isomorphic to $(W,\omega')$.
\end{defn}

Note that in this case $\jj^{\perp}$ is an ideal on $\ggo$, but this is not true in general for every  isotropic ideal on $\ggo$.

\begin{cor} If $\ggo$ is a symplectic Lie algebra isomorphic to either $\rr\hh_3$, or $\rr\rr_{3,0}$ or $\rr\rr_{3,-1}$ or $\rr\rr_{3,0}'$ or $\nn_4$ or $\rr_{4,0}$,  then $\ggo$ is a symplectic double extension of $\RR^2$.
\end{cor}
\begin{proof}  Since these symplectic Lie algebras have non trivial center, the assertion follows from the previous theorem.
\end{proof}

Let $\hh$ be a central one dimensional subalgebra of a symplectic Lie algebra $\ggo$, then the following exact sequences of Lie algebras describe central extensions of Lie algebras.
\begin{equation}\label{e1}
0 \longrightarrow \hh \longrightarrow \hh^{\perp} \longrightarrow \hh^{\perp}/\hh \longrightarrow 0
\end{equation}
\begin{equation}\label{e2}
0 \longrightarrow \hh^{\perp} \longrightarrow \ggo \longrightarrow \ggo/ \hh^{\perp}\longrightarrow 0
\end{equation}
\begin{equation}\label{e3}
0 \longrightarrow \hh^{\perp}/\hh \longrightarrow \ggo/\hh \longrightarrow \ggo/ \hh^{\perp}\longrightarrow 0
\end{equation}
\begin{equation}\label{e4}
0 \longrightarrow \hh \longrightarrow \ggo \longrightarrow \ggo/ \hh\longrightarrow 0
\end{equation}
Since $\dim \hh =1$ then (\ref{e3}) and (\ref{e4}) describe semidirect products of Lie algebras.

For the symplectic Lie algebra $\rr_{3,0}'$, let $\hh$ be the central ideal spanned by e$_4$, then $\hh^{\perp}=<e_2,e_3,e_4>$. The first exact sequence (\ref{e1}) describe a trivial central extension, that is, $\hh^{\perp}$ is the extension of $\RR^2=<e_2,e_3>$ by $\RR$ defined by the zero class in $Z^2(\hh^{\perp}/\hh,\RR)$. The second exact sequence (\ref{e2}) describe a semidirect product of Lie algebras, by the action of $\RR$ on $\hh^{\perp}$ via ade$_1$.

In order to give a model for the symplectic Lie algebras $\rr_{4,0}'$ and $\dd_{4,\delta}'$ for $\delta >0$,  we make use of generalized symplectic Lie algebras. This construction  was introduced by  Dardi\'e and Medina in \cite{DM1}. A generalized  symplectic Lie algebra $\aa$ admits a decomposition $V^{\ast} \oplus B \oplus V$ as vector spaces, where $(B, \omega')$ is a symplectic Lie algebra, $V^{\ast}$ is the dual vector space of $V$, which  is a the underlying Lie algebra of a left symmetric algebra and  such that  there exists a representation $\Gamma : V \to Der(B)$,  a cocycle $\varphi \in Z^2_{S.G.}(V, V^{\ast})$ and a symmetric bilinear form $f: V \times V \to B$, satisfying some extra conditions. These conditions assert that $\aa$ is a symplectic Lie algebra endowed with the symplectic structure $\omega' + \omega_0$ where $\omega_0$ is the canonical symplectic form on $V^{\ast} \oplus V$ (see Thm 2.3 in \cite{DM1}). We remite to this paper for more explanation. We will apply the following result of Dardi\'e and Medina \cite{DM1} to characterize some four dimensional  symplectic Lie algebras.

\begin{thm} Let $(\ggo,\omega)$ be a symplectic Lie algebra and let $\hh$ be a isotropic ideal of $\ggo$. If $\hh^{\perp}$ is a ideal of $\ggo$ and if the following exact sequence of Lie algebras
\begin{equation}\label{e5}
0 \longrightarrow \hh^{\perp}/\hh \longrightarrow \ggo/\hh \longrightarrow \ggo/ \hh^{\perp}\longrightarrow 0
\end{equation}
splits, then $\ggo$ is a generalized symplectic double extensions of the reduced symplectic Lie algebra $\hh^{\perp}/\hh$ by the Lie algebra $\ggo/ \hh^{\perp}$.\end{thm}
Thus the following exact sequences characterize the generalized symplectic double extension:
$$0  \longrightarrow  \hh \longrightarrow  \hh^{\perp} \longrightarrow \hh^{\perp}/\hh \longrightarrow 0$$
$$0  \longrightarrow  \hh^{\perp} \longrightarrow  \ggo \longrightarrow \ggo/\hh^{\perp} \longrightarrow 0$$
$$0  \longrightarrow  \hh \longrightarrow \ggo \longrightarrow \ggo/\hh \longrightarrow 0$$
$$0  \longrightarrow  \hh^{\perp}/\hh \longrightarrow \ggo/\hh \longrightarrow \ggo/\hh^{\perp} \longrightarrow 0$$

\begin{remark} According to the previous theorem,  Lie algebras having a lagrangian ideal are generalized symplectic double extensions. In fact, one can  easily check that the exact sequence of Lie algebras:
$$
0 \longrightarrow \hh^{\perp}/\hh \longrightarrow \ggo/\hh \longrightarrow \ggo/ \hh^{\perp}\longrightarrow 0
$$
splits when $\hh=\hh^{\perp}$ is an ideal.
\end{remark}

\begin{prop} Let $\ggo$ be a symplectic  Lie algebra isomorphic either to $\rr_{4,0,\delta}'$ or $ \dd_{4,\delta}'$ $\delta>0$, then $\ggo$  is a generalized symplectic double extension of $\RR^2$ by $\RR$.
\end{prop}

\begin{proof} In the cases $\rr_{4,0,\delta}'$, $ \dd_{4,\lambda}'$, let $\hh$ be the isotropic ideal generated by $e_1$ in the first case and by $e_3$ in the second one. Then the  orthogonal subspace $\hh^{\perp}$ is the ideal  $\hh^{\perp}=<e_1,e_2,e_3>$  in both cases. It is easy to see  that $\hh^{\perp}/\hh$ is  abelian and two dimensional. Thus
one has the exact sequences of Lie algebras:
\begin{equation}\label{2} 0 \to \hh \to  \hh^{\perp} \to  \hh^{\perp} / \hh  \to 0\end{equation}
\begin{equation}\label{3} 0 \to \hh^{\perp} \to \ggo \to \ggo / \hh^{\perp} \to 0\end{equation}
where $\bb:=\hh^{\perp}/\hh$ is a abelian two dimensional Lie algebra endowed with a symplectic structure induced from $\ggo$. The first sequence (\ref{2}) is defined by the cohomology class of a 2-cocycle $\varphi \in Z^2(\bb,\RR)$. If the commutator of $\ggo$ is abelian then  $\varphi$ is trivial and thus $\hh^{\perp}$ is a direct product of $\hh$ and $\bb$. If $\ggo$ is isomorphic to $\dd_{4,\lambda}'$, then $\varphi$ is no trivial. Thus for $\bb$ the abelian two dimensional Lie algebra generated by $e_1$ and $e_2$, the Lie bracket on $\hh^{\perp}=\bb \oplus \RR e_3$ is given by $[e_1, e_2]= \varphi(e_1, e_2) e_3$.

The second sequence (\ref{3}) splits (see \cite{A-B-D-O}) and so we prove that $\ggo$ is a a generalized symplectic double extension.
\end{proof}

\

\section{Some generalizations: Obstructions}\label{last}

Motivated by the results in the four dimensional case, we generalize some results. In particular propositions in this section  should be compared with those in previous sections which describe the structure of four dimensional Lie algebras admitting  symplectic structures. The antiderivation operator  in  $\Lambda(\ggo^{\ast})$ will be denoted $d$ throughout this section.

\medspace

Recall that the Heisenberg Lie algebra of dimension 2n+1, denoted $\hh_{2n+1}$, is generated by elements $e_i$ i=1, $\hdots$ 2n+1, with the relations
$[e_i,e_{i+1}]= e_{2n+1}$ i=2k+1, k=0, $\hdots$ , n-1. Using this Lie bracket relations the following lemma follows.

\begin{lem} Let  $D$ be a derivation of the Heisenberg Lie algebra $h_{2n+1}$. Then, in the basis $e_1, \hdots , e_{2n+1}$, the matrix of $D$ is
$$
\left( \begin{matrix}
A & \ast \\
0 & \lambda
\end{matrix}
\right)
\quad A =(a_{ij}) \in gl_n(\RR),\, a_{i+1,i+1} = \lambda -a_{i,i},\, i=2k+1,\, k=0, \hdots n-1.
$$
\end{lem}

It follows that $\lambda = tr D/n$.

In particular the unimodular extensions of $\hh_{2n+1}$ are those such that $\lambda = 0$. In fact

tr(D) =$ \lambda +  \sum a_{s,s} = \lambda + \sum_{r \, even} a_{r,r} + \sum_{t\, odd} a_{t,t}=\lambda + \sum_{r \, odd} (a_{r,r} +  \lambda - a_{r,r})= (n+1) \lambda$,

\noindent and this implies the assertion.

Observe that changing $D$ by $\tilde{D}:=D-\sum_{i=1}^{2n} a_{i, 2n+1}e_i$ we get a new derivation of $\hh_{2n+1}$ of the form,
\begin{equation}
\left( \begin{matrix}
A & 0 \\
0 & \lambda
\end{matrix}
\right)
\label{A} \quad A =(a_{ij}) \in gl_n(\RR),\, a_{i+1,i+1} = \lambda -a_{i,i},\, i=2k+1,\, k=0, \hdots n-1.
\end{equation}
The extended Lie algebras resulting  as semidirect products of the Heisenberg Lie algebra by $\RR$, using $D$ or $\tilde{D}$,  are isomorphic.

\begin{prop} Let $\ggo= \RR e_0 \ltimes \hh_{2n+1}$ be a unimodular extension of the Heisenberg Lie algebra $\hh_{2n+1}$ such that, $A$ as in (\ref{A}), belongs to $\GL_n(\RR)$. Then $\ggo$ does never admit a symplectic structure.
\end{prop}
\begin{proof}
It holds $de^0=0$ and if $A\in \GL_n$ then $\{z_j = Ae_j\}$ j=1, $\hdots$ , 2n, is a basis of the subspace spanned by $e_1,\hdots , e_{2n}$. Then we have
$$
de^{2n+1} = -\sum_{i\, even} e^i \wedge e^{i+1},\quad de^0 =0 \quad dz^j = -e^0 \wedge e^j, \, j=1,\hdots , 2n$$
Thus $$
d(e^0 \wedge e^{2n+1}) \ne 0, \quad\mbox{and } \quad d(z^j \wedge e^{2n+1})  \ne 0
$$
So if $\omega$ is a closed two-form then $\omega \in \Lambda^2(W^{\ast})$ where $W$ is the subspace of $\ggo$ generated by $e_i$ i = 1,..., 2n and this implies $\omega^{n+1} =0$.
\end{proof}

\begin{remark} Let $\ggo$ be a Lie algebra as in the previous proposition. Then  $H_{2n+2}(\ggo) \ne 0$. In fact let $\ad_{e_0}=D$. Then $\Lambda^{2n+2}(\ggo)$ is generated by $e_0 \wedge e_1 \wedge \hdots \wedge e_{2n+1}$ and let $\partial_i$ the corresponding coboundary operator at level i. Thus at the level 2n+2, one has:
$$\begin{array}{rcl}
\partial_{2n+2} (e_0 \wedge \hdots \wedge e_{2n+2}) &  = & \sum_{i= 1}^{2n} (-1)^i [e_0, e_i] \wedge e_1 \wedge \hdots \wedge \mbox{\^{e}}_i \wedge \hdots e_{2n+1} \\
& =  & \sum_{i \, even } (-1)^i a_{ii} e_i \wedge e_1 \wedge \hdots \wedge \mbox{\^{e}}_i \wedge \hdots e_{2n+1} + \\
& & + \sum_{j \, odd } (-1)^j a_{jj} e_j \wedge e_1 \wedge \hdots \wedge \mbox{\^{e}}_j \wedge \hdots e_{2n+1}\\
& = & 0
\end{array}
$$
where $a_{ii}$ is the entry ii of the matrix $A$ of the previous Lemma, and satisfies the corresponding condition for unimodular extensions. At the level 2n+1 it holds:
$$\begin{array}{rcl}
\partial_{2n+1} (e_1 \wedge \hdots e_{2n+1}) &  = & \sum_{i<j} (-1)^{i+j+1} [e_i, e_j] \wedge  \hdots \wedge \mbox{\^{e}}_i \wedge \hdots \wedge \mbox{\^{e}}_j \wedge \hdots \wedge e_{2n+1} \\
& = & 0
\end{array}
$$
Thus $H_{2n}(\ggo) \ne 0$. However $H^2(\ggo)$ does not necessarly vanishes. In fact, take $D$ the derivation of $\hh_{5}$ given by $De_1= e_1$, $De_e =- e_2$ $De_3=-e_3$ and $De_4 = e_4$. Then $d(e^1 \wedge e^3) =0$ and so $H^2(\ggo) \ne 0$.
\end{remark}

It is known that trivial extensions of the Heisenberg Lie algebras $\hh_{2n+1}$ are symplectic if and only if n=1. Here we give a proof of this fact.

\

\begin{prop} Let $\ggo$ be a trivial extension of the Heisenberg Lie algebra $\hh_{2n+1}$. Then $\ggo$ is symplectic if and only if n=1.
\end{prop}
\begin{proof}
One has $de^i=0$ for i=0, $\hdots$, 2n and $de^{2n+1} = -\sum_{i \, odd} e^i \wedge e^{i+1}$. Thus $d(e ^i \wedge e^j) =0$ for i=0, $\hdots$, 2n and $d(e^{2n+1} \wedge e^j) = -\sum_{i\, odd, i, i+1\ne j} e^i \wedge e^{i+1} \wedge e^j\ne 0$  for j$\ne$ 2n+1. Thus if $\omega $ is closed then $\omega$ belongs to $\Lambda^2(W^{\ast})$ for $W=\la e_0 , \hdots, e_{2n}\ra$ and that implies $\omega^{n+1} =0$ for n $\ge 2$.
\end{proof}

\begin{prop} Let $\ggo$ be a semidirect product of $\RR e_0$ and the 2n-1 dimensional abelian ideal such that $\ad_{e_0}$ diagonalizes and the eigenvalues of $\ad_{e_0}$ satisfy $\lambda_i + \lambda_j \ne 0$ for all i,j. Then $\ggo$ cannot be equipped with a symplectic structure.
\end{prop}
\begin{proof}
 Since  $de_i=-\lambda_i e^0 \wedge e^i$, then $d(e^i \wedge e^j) = -(\lambda_i + \lambda_j) e^0 \wedge e^i \wedge e^j) \ne 0$. Thus if $\omega$ is closed two-form, then $\omega = \sum \alpha_i de_i$ and so $\omega^n =0$.
\end{proof}

\begin{prop} Let $\ggo$ be a semidirect product of $\RR e_0$ and the 2n-1 dimensional abelian ideal such that $\ad_{e_0}e_i = e_{i-1}$,  $2 \le i \le 2n-1$. Then $\ggo$ is symplectic if and only if n=2.
\end{prop}
\begin{proof}
If n=2 then (\ref{simplecticas}) proves that the corresponding Lie algebra $\ggo\simeq \nn_4$ can be symplectic.  Assume that $n\ge 3$. Since  $de_i=- e^0 \wedge e^{i+1}$ and $de^0= 0 = de^{2n-1}$, then $d(e^{2n-2} \wedge e^{2n-1})=0= d(e^0 \wedge e^{i})$ for $i =1, \hdots 2n-1$, but $d(e^i \wedge e^j)\ne 0$ for $1\le i < j \le 2n-1$. Thus any closed two-form $\omega$ has the form  $\omega = be^{2n-2} \wedge e^{2n-1} + \sum \alpha_i de^0 \wedge e^i$ and so $\omega^3 =0$, which implies $\omega$ cannot be symplectic.
\end{proof}

\

{\it Acknowledgements.} The author was partially supported
by CONICET and SECYT-UNC (Argentina).

The author thanks I. Dotti  for general supervision, M. Fern\'andez for her comments and for suggesting the subject of this article and L. Cagliero for very useful discussions and comments.

\section{Appendix I: Automorphisms of symplectic Lie algebras}

We compute the automorphism of symplectic four dimensional Lie algebras according to the list obtained in Proposition (\ref{simplecticas}).
We identify a automorphism $\sigma$ with its matrix representation in  the ordered basis $\{e_1,e_2,e_3,e_4\}$ as in Prop. (\ref{clases}). In all cases the matrix must be non singular.

\medspace

\begin{tabular}{lclc}
$\rr\hh_3:$ & ${\left( \begin{matrix} a_{11} & a_{12} & 0 & 0 \\
a_{21} & a_{22} & 0 & 0 \\
a_{31} & a_{32} & a_{33} & a_{34} \\
a_{41} & a_{42} & 0 & a_{44} \\
\end{matrix}\right)}$ &
$\rr \rr_{3,0}:$ & ${\left( \begin{matrix} 1 & 0 & 0 & 0 \\
a_{21} & a_{22} & 0 & 0 \\
a_{31} & 0 & a_{33} & a_{34} \\
a_{41} & 0 & a_{43} & a_{44} \\
\end{matrix}\right)}$ \\
& {$a_{33}=a_{11} a_{22} - a_{12} a_{21}$} & & \\ \\
$\rr\rr_{3,-1}:$ & ${\left( \begin{matrix} 1 & 0 & 0 & 0 \\
a_{21} & a_{22} & 0 & 0 \\
a_{31} & 0 & a_{33} & 0 \\
a_{41} & 0 & 0 & a_{44} \\
\end{matrix}\right)}$ &
or & ${\left( \begin{matrix} -1 & 0 & 0 & 0 \\
a_{21} & 0 & a_{23} & 0 \\
a_{31} & a_{32} & 0 & 0 \\
a_{41} & 0 & 0 & a_{44} \\
\end{matrix}\right)}$ \\
\\
\end{tabular}

\pagebreak

\begin{tabular}{lclc}
$\rr\rr_{3,0}':$ & ${\left( \begin{matrix} 1 & 0 & 0 & 0 \\
a_{21} & a_{22} & a_{23} & 0 \\
a_{31} & -a_{23} & a_{22} & 0 \\
a_{41} & 0 & 0 & a_{44} \\
\end{matrix}\right)}$ &
or & ${\left( \begin{matrix} -1 & 0 & 0 & 0 \\
a_{21} & a_{22} & a_{23} & 0 \\
a_{31} & a_{23} & -a_{22} & 0 \\
a_{41} & 0 & 0 & a_{44} \\
\end{matrix}\right)}$ \\ \\
$\rr_2 \rr_{2}:$ & ${\left( \begin{matrix} 1 & 0 & 0 & 0\\
a_{21} & a_{22} & 0 & 0 \\
0 & 0 & 1 & 0  \\
0 & 0 & a_{43} & a_{44} \\
\end{matrix}\right)}$ &
or &  ${\left( \begin{matrix} 0 & 0 & 1 & 0\\
0 & 0 & a_{23} & a_{24} \\
1 & 0 & 0 & 0 \\
a_{41} & a_{42}  & 0 & 0 \\
\end{matrix}\right)}$\\ \\
$\rr_{2}':$ & ${\left( \begin{matrix} 1 & 0 & 0 & 0\\
0 & 1 & 0 & 0 \\
a_{31} & a_{32} & a_{33} & a_{34} \\
a_{41} & a_{42} & -a_{34}  & a_{33} \\
\end{matrix}\right)}$ &
or &  ${\left( \begin{matrix} 1 & 0 & 0 & 0\\
0 & -1 & a_{23} & a_{24} \\
a_{31} & a_{32} & a_{33} & a_{34} \\
a_{41} & a_{42}  & a_{34} & -a_{33} \\
\end{matrix}\right)}$\\ \\
$\nn_{4}:$ & ${\left( \begin{matrix} a_{11} & 0 & 0 & a_{14}\\
a_{21} & a_{11}a_{44} & 0 & a_{24} \\
a_{31} & a_{32} & a_{33} & a_{34} \\
0 & 0 & 0 & a_{44} \\
\end{matrix}\right)}$ &
$\rr_{4,\varepsilon}:$ & ${\left( \begin{matrix} a_{11} & 0 & 0 & a_{14} \\
0  & a_{22} & a_{23} & a_{24} \\
0 & 0 & a_{22} & a_{34} \\
0 & 0 & 0 & 1 \\
\end{matrix}\right)}$ \\
& $a_{32}= a_{44}a_{21},\, a_{33}= a_{44}^2 a_{11}$ & & $\varepsilon = 0,-1 $\\ \\
$\rr_{4,-1,-1}$ & ${\left( \begin{matrix} a_{11} & 0 & 0 & a_{14} \\
0  & a_{22} & a_{23} & a_{24} \\
0 & a_{32} & a_{33} & a_{34} \\
0 & 0 & 0 & 1 \\
\end{matrix}\right)}$ &
$\rr_{4,\mu,\nu}:$ & ${\left( \begin{matrix} a_{11} & 0 & 0 & a_{14} \\
0  & a_{22} & 0 & a_{24} \\
0 & 0 & a_{33} & a_{34} \\
0 & 0 & 0 & 1 \\
\end{matrix}\right)}$ \\
&  & & {$(\mu, \nu) = (-1,\beta), (\alpha, -\alpha)$} \\ \\
$\rr_{4,0,\delta}':$ & ${\left( \begin{matrix} a_{11} & 0 & 0 & a_{14} \\
0  & a_{22} & a_{23} & a_{24} \\
0 & -a_{23} & a_{22} & a_{34} \\
0 & 0 & 0 & 1 \\
\end{matrix}\right)}$ &
$\hh_{4}:$ & ${\left( \begin{matrix} a_{11} & a_{12} & 0 & a_{14} \\
0  & a_{22} & 0 & a_{24} \\
2a_{11} a_{24} & a_{32} & a_{11}^2 & a_{34} \\
0 & 0 & 0 & 1 \\
\end{matrix}\right)}$ \\
& $\delta>0$ & & $a_{32}= 2(a_{24}(2a_{11}+a_{12}) -a_{14}a_{11})$\\ \\
$\dd_{4,1/2}:$ & ${\left( \begin{matrix} a_{11} & a_{12} & 0 & a_{14} \\
a_{21}  & a_{22} & 0 & a_{24} \\
a_{31} & a_{32} & a_{33} & a_{34} \\
0 & 0 & 0 & 1 \\
\end{matrix}\right)}$ &
$\dd_{4,\lambda}:$ & ${\left( \begin{matrix} a_{11} & 0 & 0 & a_{14} \\
0  & a_{22} & 0 & a_{24} \\
a_{31} & a_{32} & a_{11}a_{22} & a_{34} \\
0 & 0 & 0 & 1 \\
\end{matrix}\right)}$ \\ \\
& $a_{31}=2(a_{11}a_{24}-a_{14} a_{21}),$ & $\lambda > 1/2$& $a_{31}=\frac{a_{11}a_{24}}{1-\lambda},\,\mbox{ iff } \lambda \neq 1$\\
& $a_{32}=2(a_{12}a_{24} -a_{22}a_{14}),$ & & $a_{32}= -\frac{a_{14}a_{22}}{\lambda}$\\
& $a_{33}= a_{12}a_{21} -a_{22} a_{11} $ & &
\end{tabular}

$$\dd_{4,\delta}':\quad {\left( \begin{matrix} a_{11} & a_{12} & 0 & a_{14} \\
-a_{12}  & a_{11} & 0 & a_{24} \\
a_{31} & a_{32} & a_{11}^2 + a_{22}^2 & a_{34} \\
0 & 0 & 0 & 1 \\
\end{matrix}\right)}$$
 with
$$\left( \begin{matrix} a_{31}\\
a_{32} \end{matrix}\right) = \frac 1{(\delta/2)^2+1}
\left( \begin{matrix} \delta/2 & 1 \\
-1 & \delta/2
\end{matrix}\right)  \left( \begin{matrix} a_{12} & a_{11} \\
-a_{11} & a_{12}
\end{matrix}\right)  \left( \begin{matrix} a_{14} \\
a_{24}
\end{matrix}\right) $$

\section{Appendix II: Cohomology} To conclude and for further use we compute the cohomology of these Lie algebras.

\begin{prop} \label{coho} The cohomology over $\RR$ of any four dimensional solvable real Lie algebra is presented in the following table:

\begin{center} 
\small{
\begin{tabular}{|c|c|c|c|}\hline 
{{\rm Case}} & {{$H^1(\ggo)$}} & {${H^2(\ggo)}$} & {$H^3(\ggo)$}\\ \hline 
{$\rr \hh_3$} & {$[e^1][e^2][e^4]$} & {$[e^1 \wedge e^3][e^1 \wedge e^4]$} & {$[e^1 \wedge e^2 \wedge e^3] [e^1 \wedge e^3 \wedge e^4] $} \\ 
 & & {$[e^2 \wedge e^3][e^2 \wedge e^4]$} & {$[e^2 \wedge e^3 \wedge e^4]$} \\ \hline
{$\rr\rr_{3}$} & {$[e^1][e^4]$} & {$[e^1 \wedge e^4]$} & {$0$} \\
\hline
{$\rr \rr_{3,0}$} & {$[e^1][e^4]$} & {$[e^1 \wedge e^4]$} & {$[e^1 \wedge e^2 \wedge e^3]$} \\
\hline
{$\rr \rr_{3,-1}$} & {$[e^1][e^4]$} & {$[e^1 \wedge e^4][e^2 \wedge e^3]$} & {$[e^1 \wedge e^2 \wedge e^3][e^2 \wedge e^3 \wedge e^4]$} \\
\hline
{$\rr \rr_{3,\lambda},\,\lambda \ne 0,-1$} & {$[e^1][e^4]$} & {$[e^1 \wedge e^4]$} & {$0$} \\
\hline
{$\rr \rr'_{3,0}$} & {$[e^1][e^4]$} & {$[e^1 \wedge e^4][e^2 \wedge e^3]$} & {$[e^1 \wedge e^2 \wedge e^3][e^2\wedge e^3\wedge e^4] $} \\ 
\hline
{$\rr \rr'_{3,\gamma}$, $\gamma \ne 0$} & {$[e^1][e^4]$} & {$[e^1 \wedge e^4]$} & {$0$} \\ 
\hline
{$\rr_2 \rr_2$} & {$[e^1][e^3]$} & {$[e^1 \wedge e^3]$} & {$0$} \\ 
\hline
{$\rr_2'$} & {$[e^1][e^2]$} & {$[e^1 \wedge e^2]$} & {$0$} \\ 
\hline
{$\nn_4$} & {$[e^1][e^4]$} & {$[e^1 \wedge e^2][e^3 \wedge e^4]$} & {$[e^1 \wedge e^2 \wedge e^3] [e^2 \wedge e^3 \wedge e^4] $} \\
\hline
{$\rr_4$} & {$[e^4]$} & {$0$} & {$0$} \\
\hline
{$\rr_{4,0}$} & {$[e^3][e^4]$} & {$[e^2 \wedge e^3][e^2 \wedge e^4]$} & {$ [e^2 \wedge e^3 \wedge e^4]$} \\
\hline
{$\rr_{4,-1}$} & {$[e^4]$} & {$[e^2 \wedge e^4]$} & {$0$} \\ 
\hline
{$\rr_{4,-1/2}$} & {$[e^4]$} & {$0$} & {$[e^1 \wedge e^2 \wedge e^3]$} \\ 
\hline
{$\rr_{4,\mu},\, \mu\ne -1,-1/2,0$} & {$[e^4]$} & {$0$} & {$0$} \\ 
\hline
{$\rr_{4,-1, -1}$} & {$[e^4]$} & {$[e^1 \wedge e^3]$} & {$[e^1 \wedge e^2 \wedge e^4] [e^1\wedge e^3\wedge e^4] $} \\
\hline 
{$\rr_{4,-1, \beta}$} {$\beta \ne -1$}  
& {$[e^4]$} & {$[e^1 \wedge e^2]$} & {$[e^1 \wedge e^2 \wedge e^4]  $} \\
\hline 
 

{$\rr_{4,\alpha, -\alpha}$} {$\alpha \ne -1,0$}
& {$[e^4]$} & {$[e^2 \wedge e^3]$} & {$[e^2 \wedge e^3 \wedge e^4]  $} \\\hline
{$\rr_{4,\alpha, \beta}\,$}& & {} & \\
{$\alpha \ne -1,0,-\beta $},{$\alpha - \beta = -1$}
 & {$[e^4]$} & {$0$} & {$[e^1\wedge e^2\wedge e^3]  $} \\ \hline
{$\rr_{4,\alpha, \beta}, \,\,\alpha, \beta$}& & {} & \\
{\rm not as above }  {$\alpha - \beta \ne -1$}& {$[e^1]$} & {$0$} & {$0$} \\ \hline
{$\rr'_{4,0, \delta}$} {$\delta \ne 0$}& {$[e^4]$} & {$[e^2 \wedge e^3]$} & {$ [e^2 \wedge e^3 \wedge e^4]$} \\ 
\hline 
{$\rr'_{4,-1/2, \delta}$} {$\delta \ne 0$}& {$[e^4]$} & {$0$} & {$ [e^1 \wedge e^2 \wedge e^3] $} \\ 
\hline
{$\rr'_{4,\gamma, \delta}$} {$\gamma \ne -1/2, 0, \delta \ne 0$}& {$[e^4]$} & {$0$} & {$0$} \\
\hline
{$\dd_4$} & {$[e^4]$} & {$0$} & {$[e^1 \wedge e^2 \wedge e^3]$} \\
\hline
{$\dd_{4,1}$} & {$[e^2][e^4]$} & {$ [e^2 \wedge e^4]$} & {$0$} \\
\hline
{$\dd_{4,2}$} & {$[e^4]$} & {$[e^2 \wedge e^3]$} & {$[e^2 \wedge e^3 \wedge e^4] $} \\ 
\hline
{$\dd_{4,\lambda}$}  {$ \lambda  \ne 1,2\quad $}& {$ [e^4]$} & {$0$} & {$0$} \\
\hline
{$\dd'_{4,0}$} & {$[e^4]$} & {$0$} & {$[e^1 \wedge e^2 \wedge e^3][e^1\wedge e^2 \wedge e^4]$} \\
\hline
{$\dd'_{4,\delta}$}{$\delta \ne 0$} & {$[e^4]$} & {$0$} & {$0$} \\ 
\hline
{$\hh_4$} & {$[e^4]$} & {$0$} & {$0$} \\
\hline
\end{tabular} }
\end{center}
\begin{center} 
{\rm Table} \ref{coho} 
\end{center} 
\end{prop}
\begin{proof} The cohomology can be obtained parallel to the computations made to get the Table (\ref{simplecticas}). Continuing with the case $\rr_2'$, (worked out in the previous Proposition (\ref{simplecticas})), in the proof of this Proposition we can see that $H^1(\ggo) =\{ [e^1], [e^2]\}$. From the computations at the next level it is possible to prove that $\theta \in H^2(\ggo)$ if and only if $\theta$ belongs to the class $[e^1 \wedge e^2]$. Since $e^1 \wedge e^2 \wedge e^3, e^1 \wedge e^2 \wedge e^4$ and $e^1 \wedge e^3 \wedge e^4$ are in the image of $d:\Lambda^{2}(\ggo) \to \Lambda^{3}(\ggo)$, to get $H^3(\ggo)$ we need to compute extra only the following: $d(e^2 \wedge e^3 \wedge e^4)= 2 e^1 \wedge e^2 \wedge e^3 \wedge e^4$ and thus we get the results of the Table for this case. The other cases can be handled in a similar way to complete the proof of the Table (\ref{coho}).
\end{proof}

\end{document}